\documentclass[11pt,letterpaper]{article}  %[journal,onecolumn,10pt]{IEEEtran}
\usepackage{epsfig}
\usepackage{amsfonts}
\usepackage{amssymb}
\usepackage{amsmath}
\usepackage{float}
\usepackage{xcolor}
\newcommand{\ba}{\begin{array}}
\newcommand{\ea}{\end{array}}
\newcommand{\bal}{\begin{align}}
\newcommand{\eal}{\end{align}}
\newcommand{\be}{\begin{equation}}
\newcommand{\beqn}{\begin{equation}}
\newcommand{\eeqn}{\end{equation}}
\newcommand{\bea}{\begin{eqnarray}}
\newcommand{\eea}{\end{eqnarray}}
\newcommand{\benum}{\begin{enumerate}}
\newcommand{\eenum}{\end{enumerate}}
\newcommand{\bi}{\begin{itemize}}
\newcommand{\ei}{\end{itemize}}
\newcommand{\bean}{\begin{eqnarray*}} %non numbered eqns.
\newcommand{\eean}{\end{eqnarray*}}   %non numbered eqns.
\def\arg{\mbox{\rm arg}}

\def\qed{\hfill$\Box$}
\def\bS{{\bf S}}
\newtheorem{theorem}{Theorem}
\newtheorem{defn}{Definition}
\newtheorem{lemma}{Lemma}

\newtheorem{example}{Example}
\newtheorem{remark}{Remark}
\def\proof{{\it Proof.\ }}
\begin{document}

\vspace{-0.3cm}

\title{\LARGE\bf Mean Field Stochastic Games with Binary Action Spaces and Monotone Costs}

%\thanks{A conference  version has been  presented at the 55th IEEE CDC conference, %2016.}}

%\footnote{M. Huang's work was supported in part by Natural Sciences and Engineering Research Council of Canada %under a discovery grant and a discovery accelerator supplement grant.}}

\author{Minyi Huang\thanks{M. Huang is with the School of Mathematics and Statistics, Carleton University, Ottawa, ON
K1S 5B6, Canada (mhuang@math.carleton.ca). }
\quad and\quad  Yan Ma\thanks{Y. Ma is with the School of Mathematics and Statistics, Zhengzhou University, 450001,  Henan, China (mayan203@zzu.edu.cn). }
}

  \date{} %Feb., 2016}
\maketitle

\begin{abstract}
This paper considers  mean field games in a multi-agent Markov decision process (MDP) framework.
Each player has a continuum  state  and  binary action. By active control, a player can bring its state to a resetting point.
 All players are coupled through their cost functions. The structural property of the individual strategies is characterized in terms of threshold policies when the mean field game admits a solution.
 We further introduce a stationary  equation system of the mean field game and
 analyze uniqueness of its solution under positive externalities.
\end{abstract}

\begin{center}
\begin{minipage}{15.3cm}
\end{minipage}
\end{center}

%\begin{center}
%\begin{minipage}{15.3cm}
%{\bf AMS subject classifications:}
%60J05,  %      Discrete-time Markov processes on general state spaces
%90C40,  %      Markov and semi-Markov decision processes
%91A10,  %      Noncooperative games
%91A15  %       Stochastic games     %91A25,  %         Dynamic games
% %93E20, %     Optimal stochastic control
%\end{minipage}
%\end{center}

\section{Introduction}

Mean field game theory studies stochastic decision problems with
a large number of noncooperative players
which are individually insignificant but collectively have a significant impact on a particular player. It provides a powerful methodology for reducing complexity in the analysis and design of strategies.  With the aid of an infinite population model, one may apply consistent mean field approximations to
construct a set of decentralized strategies
for the original large but finite population model and show its $\varepsilon$-Nash equilibrium property \cite{HCM03,HCM07,HMC06}. A closely related approach is independently developed in \cite{LL07}.
 Another related solution notion in Markov decision models is the oblivious equilibrium \cite{WBR08}.
For nonlinear diffusion models
\cite{C12,HMC06,LL07},  the analysis of mean field games depends on tools of Hamilton-Jacobi-Bellman  (HJB) equations, Fokker-Planck equations, and McKean-Vlasov equations.
For further literature in the stochastic analysis setting, see \cite{CD13,KLY11}.
To address mean field interactions with an
agent possessing strong influences, mixed player models  are
studied in \cite{BFY13,H10,NH12,NC12}.  The readers are referred to \cite{BFY13,C14,GS14} for an overview on mean field game theory.

Mean field games have found applications in diverse areas such as
 power systems \cite{KM13}, large population electric vehicle recharging control \cite{MCH13,PC14}, economics and finance \cite{AJW15,CS14,LM13},
stochastic growth theory \cite{H13}, bio-inspired oscillator games \cite{YMMS12}.

This paper studies a class of mean field games in a multi-agent Markov decision process (MDP) framework.
Dynamic games within an MDP setting are a classic area pioneered by Shapley  under the name stochastic games \cite{FV97,S53}. For MDP based mean field game modeling, see \cite{AJW15,HMC04,WBR08}.
 The players in our model have continuum state spaces and binary action spaces, and have  coupling through their cost  functions. The state of each player is used to model its risk (or distress) level which has random increase if no active control is taken. The one stage cost of a player depends on its own state, the population average state and its control effort. Naturally, the cost of a player is an increasing function of its own state. The motivation of this modeling framework comes from  applications including network security investment games, and flue vaccination games \cite{BE04,JAW11,LB08,MP10};
when the cost function is  an increasing function of the population average state, it reflects positive externalities. Markov decision processes with binary action spaces also arise in control of queues
and machine replacement problems \cite{AS95,BR11}.
Our game model has connection with anonymous sequential games
\cite{JR88} which combine stochastic game modeling with a continuum of players.
However, there is a subtle difference regarding
treating individual behavior.
 In anonymous sequential games one determines the equilibrium as a joint state-action distribution of the  population and  leaves the individual
strategies unspecified \cite[Sec. 4]{JR88}, although there is an interpretation of randomized actions for players sharing a given state.
Our approach works in the other direction by explicitly specifying
the best response of an individual and  using its
closed-loop state distribution to determine the mean field.
Our  modeling starts with a finite population, and  avoids certain measurability difficulties in directly treating a continuum of random processes \cite{A04}.

A very interesting feature of our model is threshold  policies for the solution of the mean field game.
We consider a finite time horizon game, and identify conditions for the existence of  a solution to the fixed point problem. The further analysis deals with the stationary equation of the game
and addresses uniqueness under positive externalities, which is done by studying ergodicity of the closed-loop state process of an individual player. Proving uniqueness results in mean field games
is a nontrivial task, particularly when attempting to seek less restrictive
conditions. This work is perhaps the first to  establish uniqueness by the route of exploiting externalities. For this paper, in order to maintain a balance in analyzing the finite horizon problem and the stationary equation system, the existence analysis of the latter is not included and will be reported in another work.

Although mean field games provide a powerful paradigm for substantially reducing complexity in designing strategies, except for the linear-quadratic (LQ) cases \cite{HCM03,LZ08,TZB11,WZ12}  allowing simple computations,
  strategies in general nonlinear systems are often only implicitly determined, rarely taking simple forms.  Their  numerical solutions lead to high computational load.  One of the objectives in this paper is to develop a modeling framework to obtain relatively simple solution structures.

This paper is an English version of
\cite{HM16Chen}. All assumptions and results in Sections 2-6 and Appendices A and B of both papers are the same, but Appendices C and D have  been rewritten.

The organization of the paper is as follows.
The Markov decision process framework  is introduced in
Section \ref{sec:mean}. Section \ref{sec:limit} solves the best
response as a threshold policy.
Section \ref{sec:solution} shows an $\epsilon$-Nash equilibrium
 property. The existence
of a solution to the mean field equation system is analyzed in
Section \ref{sec:ex}.
Section \ref{sec:se} introduces the stationary equation system and analyzes uniqueness of
its solution.
 Section \ref{sec:con} concludes the paper.

%\newpage

\section{The  Markov Decision Process Model}
\label{sec:mean}

\subsection{The system dynamics}

The system consists of $N$  players
denoted by ${\cal A}_i$, $1\le i\le N$.
At time $t\in \mathbb{Z}_+=\{0, 1,2, \ldots\}$, the state of ${\cal A}_i$ is denoted by $x_t^i$, and its action by $a_t^i$.

For simplicity, we consider a population of homogeneous players.
Each player has a state space ${\bf S}=[0,1]$.
 A value of ${\bf S}$ may be interpreted as a risk (or ``distress") level.
All players have the same action space  ${\bf A}=\{a_0, a_1\}$.
A player can either
 do nothing (action $a_0$) or
 make an active effort (action $a_1$).
  For an interval $I$, let ${\cal B}(I)$ denote the Borel $\sigma$-algebra of $I$.

The state of each player evolves as a controlled Markov  process
 which is affected only by its own action.
For $t\ge 0$ and $x\in \bS$,  the state has  transition kernel specified  by
\begin{align}
&P(x_{t+1}^i\in B|x_t^i =x, a_t^i=a_0)= Q_0(B|x), \label{xa0}\\
&P(x_{t+1}^i=0|x_t^i =x, a_t^i=a_1)=1 , \label{xa1}
\end{align}
where $Q_0(\cdot|x)$ is a stochastic kernel defined for
$B\in {\cal B}(\bS) $
and $Q_0([0,x)|x)=0$.
The structure of $Q_0$ indicates that given $x_t^i=x$, $a_t^i=a_0$, the state
has a transition into $[x, 1]$.
 In other words, the state of the player deteriorates if no active control is taken.
We call $a_t^i=a_1$ and $0\in \bS$   a resetting action and a resetting point, respectively.

The vector process $(x_t^1, \ldots x_t^N)$  constitutes a controlled Markov process in higher dimension with its transition kernel determined as a product measure of the form
$$
P(x_t^i\in B_i,i=1, \ldots, N|x_t^i=x^{[i]}, a_t^i=a^{[i]}, i=1, \ldots, N) =\prod_{i=1}^N P(x_t^i\in B_i|x_t^i=x^{[i]}, a_t^i=a^{[i]}),
$$
where $B_i\in {\cal B}(\bS)$, $x^{[i]}\in \bS$ and
 $a^{[i]}\in {\bf A}$. This product measure implies independent transitions of the $N$ controlled Markov  processes.

\subsection{The individual costs}
Define the population average state
 $x^{(N)}_t= \frac{1}{N} \sum_{i=1}^N x_t^i$.
For ${\cal A}_i$, the one stage cost is given by
$$
c(x_t^i, x^{(N)}_t, a_t^i)= R(x_t^i, x^{(N)}_t) +\gamma 1_{\{a_t^i=a_1\}},
$$
where $\gamma>0$ and  $\gamma 1_{\{a_t^i=a_1\}}$ is the effort cost. The function $R\ge 0$  is defined on   $\bS\times\bS$ and models the risk-related cost.
For  $0<T<\infty$ and discount factor $\rho \in (0, 1)$, define the cost
\begin{align}
J_i= E\sum_{t=0}^T \rho^t c(x_t^i, x^{(N)}_t, a_t^i),
\quad 1\le i\le N. \label{jin}
\end{align}

The following assumptions are introduced.

(A1) $\{x_0^i, i\ge 1\}$ are i.i.d. random variables taking values in $\bS$
and $Ex_0^i=m_0$.

%When for each $i$, $a_t^i$ depends only on $(t,x_0^i, \ldots, x_t^i)$,
%the $N$ processes $\{x_t^i, t\ge 0\}$, %$i=1, \ldots, N$, are independent.

(A2) $R(x,z)$ is a continuous function  on  $\bS\times \bS$.
  %For each fixed $x$, $R(x, \cdot)$ %is strictly increasing,
  For each fixed $z$, $R(\cdot, z)$ is  strictly increasing.

(A3) There exists a random variable $\xi$ taking values in ${\bf S}$ such that the measure
$Q_0( \cdot |x)$ is the distribution of the random variable
$x+(1-x) \xi.$ Furthermore, $P(\xi=1)<1$.

Denote the distribution function of $\xi$ by $F_\xi$.
To avoid triviality, we assume $P(\xi =1)<1$ in (A3).

We give some motivation about (A3).
For $\bS=[0,1]$, $1-x$ is the state margin from the maximal state 1.
Denote $m_t^i=1-x_t^i$. Then given  $m_t^i$ and $a_t^i=0$, the state margin decays to $m_{t+1}^i=(1-\xi_{t+1}^i)m_t^i$ where $\xi_{t+1}^i$ has the
same distribution as $\xi$. In other words,
the  state margin decays exponentially  if no active control is applied.

\begin{remark}
In fact, (A3) implies the so-called stochastic monotonicity condition
(see e.g. \cite{AS95, BR11})
for $Q_0$. We assume the specific form of $Q_0$, aiming to obtain more
refined  properties for the resulting threshold policies and sample path behavior of players.
\end{remark}

\begin{example}
%{\bf Example:}
\emph{A lab consists of $N$ networked
computers $M_i$, $1\le i\le N$, each of which is assigned to a primary user $U_i$, $1\le i\le N$, and is occasionally accessed by others for its specific resources. A computer has an unfitness state $x_t^i\in [0,1]$, which randomly degrades due to daily use and potential exposure to malwares, etc.
A user $U_i$ can take a maintenance action $a_1$ on $M_i$ by installing or updating security software, scanning and cleaning up disk, freeing up memory space, etc., to bring it to an ideal condition $x_t^i=0$.
 The one stage cost of $U_i$ is $R(x_t^i, x_t^{(N)}) +\gamma 1_{\{a_t^i=a_1\}}$ where the dependence on $x_t^{(N)}$ is due to machine sharing and potential malware spreading
 from other machines. This model is called a labmate game.}
\end{example}

\section{The Mean Field Limit Model}

 \label{sec:limit}

\subsection{The optimal control problem}

Assume (A1)-(A3) for Section \ref{sec:limit}.
A sequence $(b_s, \ldots, b_t)$, $s\le t$, is  denoted as $b_{s,t}$.
Let $x_t^i$ be  given by \eqref{xa0}-\eqref{xa1}.
Let $x^{(N)}_t$ be approximated by a deterministic value $z_t$.
Define
$$
\bar J_i(z_{0,T}, a^i_{0,T})= E\sum_{t=0}^T \rho^t c(x_t^i, z_t, a_t^i).
$$

We call $a_t^i$ a pure Markov  policy (or strategy) at $t$ if $a_t^i(x)$ is
a mapping from ${\bf S}$ to ${\bf A}$.
We say that  $a_t^i$  is a threshold policy
with parameter $r\in [0,1]$ if
$a_t^i(x)=a_1$ for $x\ge r$ and $a_t^i(x)=a_0$ for $x<r$; this gives a feedback policy.
The analysis below will identify  properties of the optimal policy.

\subsection{The dynamic programming equation}

Denote $a_{s,t}^i=(a_s^i, \ldots, a_t^i)$ for $s\le t$. Fix the sequence $z_{0,T}$, where each $z_t\in [0,1]$.
For   $0\le s\le T$ and $x\in \bS$,  define
$$
\bar J_i(s, x, z_{0,T},a_{s,T}^i)=E\left[\sum_{t=s}^T \rho^{t-s} c(x_t^i, z_t, a_t^i)\Big|x_s^i=x\right].
$$

Define
the value function
$ V(t,x)=\inf_{a_{t,T}^i} \bar J_i(t,x,z_{0,T},a^i_{t,T})$, where $a_{0,T}^i$ is from the set of all Markov policies.
The dynamic programming equation takes the form
\begin{align}
\begin{cases}
V(t, x)=\min_{a_t^i} \left[ c(x, z_t, a_t^i)+\rho E[ V(t+1, x_{t+1}^i)|x_t^i=x]\right],\\
 V(T,x)=R(x,z_T), \qquad 0\le t<T.
\end{cases} \label{dpe}
\end{align}
Write \eqref{dpe} in the equivalent form
\begin{align}
\begin{cases}
V(t,x) =  \min \Big[\rho \int_0^1 V(t+1, y) Q_0(dy|x) + R(x,z_t), \quad
 \rho V(t+1, 0) + R(x,z_t)+ \gamma\Big] , \\
 V(T,x)= R(x, z_T),\qquad 0\le t<T.
 \end{cases}  \label{dp}
\end{align}

Denote
\begin{align}
G_t(x)=\int_0^1 V(t, y) Q_0(dy|x),   \quad 0\le t\le T. \label{G}
\end{align}

\begin{lemma}\label{lemma:con}
For each $0\le t\le T$,  $V(t, x)$ is continuous on ${\bf S}$.
\end{lemma}

\proof We prove by induction.  $V(T, x)$ is  a  continuous function of $x\in \bS$.
Suppose that $V(k, x)$ is a continuous function of $x$ for $0<k\le T$.
By (A3),
\begin{align}
G_k(x) &= \int_0^1 V(k, x+(1-x)y)dF_\xi(y) %\nonumber \\
                          =  \int_0^1 V(k, (1-y)x +y) dF_\xi(y),\label{vdfxi}
\end{align}
which combined with the induction hypothesis implies that
$
\rho G_k(x)
$
is  continuous in $x$.

Note that $R(x, z_{k-1})$ is continuous in $x$.
On the other hand, if $g_1(x)$ and $g_2(x)$ are continuous on $[0,1]$,
$\min \{g_1(x), g_2(x)\}$ is a continuous function of $x$.
 It follows from \eqref{dp} that $V(k-1, x)$ is a continuous function of $x$.

By induction,  we conclude that $V(t,x)$ is continuous in $x$ for all $0\le t\le T$. \qed

\begin{lemma} \label{lemma:in}
For each $0\le t\le T$,  $V(t, x)$ is strictly increasing on $ \bS$.
\end{lemma}

\proof
For $t=T$,  $V(T, x_1)<V(T, x_2)$  whenever $x_1<x_2$.
Suppose that for $0<k\le T$,
\begin{align}
V(k, x_1)<V(k, x_2), \qquad \mbox{for}\ x_1<x_2. \label{vkle}
\end{align}
For $0\le x_1<x_2\le 1$,
 $$
 R(x_1,z_{k-1}) <  R(x_2,z_{k-1}).
$$
  By \eqref{vdfxi} and \eqref{vkle},
\begin{align*}
\rho G_k(x_1) + R(x_1,z_{k-1})
< \rho G_k(x_2) + R(x_2,z_{k-1}).
\end{align*}
For $\alpha_1<\alpha_2$ and $\beta_1<\beta_2$, we have
  %\begin{align}
$\min\{\alpha_1, \beta_1\} <\min\{\alpha_2, \beta_2\}$.  %\nonumber \label{ab}
  %\end{align}
 Taking
\begin{align*}
&\alpha_i=\rho G_k(x_i)  + R(x_i, z_{k-1}),\quad
\beta_i= \rho V(k, 0) + R(x_i,z_{k-1})+ \gamma,
\end{align*}
we obtain $V(k-1, x_1)<V(k-1, x_2).$
By induction, $V(t, x)$ is strictly increasing for all $0\le t\le T$.
  \qed

\begin{lemma} \label{lemma:G} For $0\le t\le T$,
$G_t(x)$ is continuous and  strictly increasing in $x$.
 \end{lemma}

 \proof The lemma follows from Lemmas \ref{lemma:con} and \ref{lemma:in}, \eqref{vdfxi}
 and  $ P(\xi=1)<1$ in (A3). \qed

% We determine the optimal control
% by checking  the relation between $\rho G_{t+1}( x)$ and  $\rho V(t+1, 0)+\gamma$.

\begin{lemma} \label{lemma:pa}
For $t\le T-1$,  if
\begin{align}
\rho G_{t+1}(0)<\rho V(t+1, 0) +\gamma <\rho G_{t+1}(1), \label{g01x}
\end{align}
 there exists a unique $x^*\in (0, 1)$ such that
$
\rho G_{t+1}(x^*)=  \rho V(t+1, 0)+\gamma.
$
\end{lemma}

\proof The lemma follows from Lemma \ref{lemma:G} and  the intermediate value theorem. \qed

\begin{theorem} \label{theorem:a}
 If $t=T$, define $a_T^i=a_0$.
For $t\le T-1$, define the  policy $a_t^i(x)$ by the following rule.

i) If
$\rho G_{t+1}(1) \le \rho V(t+1, 0) +\gamma,$
take $a_t^i(x)= a_0$ for all $x\in \bS$.

ii) If
$
\rho G_{t+1}(0)\ge \rho V(t+1, 0) +\gamma,
$
take $a_t^i(x)=a_1$ for all $x\in \bS$.

iii) If  \eqref{g01x} holds,
 take $a_t^i$ as a threshold policy with parameter $x^*$ given in Lemma \ref{lemma:pa}.

 Then $a_{0,T}^i$ is an optimal policy.
\end{theorem}

\proof It is easy to see that $a_T^i=a_0$ is optimal.
Consider $t\le T-1$.
By Lemma \ref{lemma:G} and \ref{lemma:pa},
we can verify that the minimum in \eqref{dp} is attained when $a_t^i$ is chosen according to i)-iii).  \qed

\section{Solution of the Mean Field Game}

\label{sec:solution}

Assume (A1)-(A3).
To obtain a solution of the mean field game,
we introduce the  equation system
\begin{align}
\begin{cases}
V(t,x) =  \min \Big[\rho \int_0^1 V(t+1, y) Q_0(dy|x) + R(x,z_t), \quad
 \rho V(t+1, 0) + R(x,z_t)+ \gamma\Big] ,\\
    \hskip 3cm 0\le t<T \\
 V(T,x)= R(x, z_T),\\
 z_t= Ex_t^i, \quad 0\le t\le T.
 \end{cases}  \label{dpz}
\end{align}
By (A1), $z_0=m_0$. We look for a solution $(\hat z_{0,T},\hat a_{0T}^i)$ for \eqref{dpz} such that
 $\{x_t^i, 0\le t\le T\}$ is generated by
$\{\hat a_t^i(x), 0\le t\le T\}$ satisfying the rule in Theorem
\ref{theorem:a} after setting $z_{0,T}=\hat z_{0,T}$.
The last equation is the standard consistency condition in mean field games.

Consider the game of $N$ players specified by \eqref{xa0}-\eqref{jin}.
Denote $a^{-i}_{0,T}= (a_{0,T}^1, \ldots, a^{i-1}_{0,T}, a_{0,T}^{i+1}, \ldots, a_{0,T}^N)$.
  Write $J_i=J_i(a^i_{0,T}, a^{-i}_{0,T})$.

For the performance estimates, we consider the perturbation of $a_t^i$ in a  strategy space ${\cal U}_t$ consisting of all pure Markov strategies depending on $(x_t^1, \ldots, x_t^N)$.

\begin{defn}
A set of strategies $\{a_{0,T}^i,1\le i\le N \}$ for the $N$ players is called an $\epsilon$-Nash equilibrium with respect to the costs $\{J_i, 1\le i\le N\}$, where $\epsilon\ge 0$, if for any $1\le i\le N$,
$$
J_i(a^i_{0,T}, a_{0,T}^{-i} ) \le J_i(b_{0,T}^i, a_{0,T}^{-i})+\epsilon,
$$
for any $b_{0,T}^i \in \prod_{t=0}^T{\cal U}_t$.
\end{defn}

\begin{theorem}
Suppose that \eqref{dpz} has a solution $(\hat z_{0,T}, \hat a_{0,T}^{i}  )$.
Then $(\hat a_{0,T}^1, \ldots, \hat a_{0,T}^N)$ is an $\epsilon$-Nash equilibrium, i.e.,
\begin{align}
J_i(\hat a_{0,T}^i, \hat a_{0,T}^{-i})-\epsilon   \le  \inf_{a_{0,T}^i} J_i( a_{0,T}^i, \hat a_{0,T}^{-i})\le J_i(\hat a_{0,T}^i,
\hat a_{0,T}^{-i}), \quad 1\le i\le N, \nonumber
\end{align}
where $ a_{0,T}^i \in \prod_{t=0}^T{\cal U}_t$ and  $\epsilon \to 0$ as $N\to \infty$.
\end{theorem}

\proof
For $(a_{0,T}^i, \hat a_{0,T}^{-i})$,
denote the corresponding states by
$x_t^i$, and $\hat x_t^j$, $j\ne i$.  We have
\begin{align}
\lim_{N\to \infty}\max_{0\le t\le T}|x_t^{(N)}- \hat z_t|=0, \quad a.s.
\label{znz}
\end{align}
where $x_t^{(N)}= \frac{1}{N} (\sum_{j\ne i} \hat x_t^j +x_t^i).$
Denote
$$
\epsilon_{1,N}=\sup_{a_{0,T}^i}|J_i( a_{0,T}^i, \hat a_{0,T}^{-i})-\bar J_i(\hat z_{0,T}, a_{0,T}^i)|.
$$
Then by \eqref{znz}, $\lim_{ N\to \infty} \epsilon_{1,N}=0$.
Furthermore,
\begin{align*}
J_i( a_{0,T}^i, \hat a_{0,T}^{-i})=&\bar J_i(\hat z_{0,T}, a_{0,T}^i) +J_i( a_{0,T}^i, \hat a_{0,T}^{-i})-\bar J_i(\hat z_{0,T}, a_{0,T}^i)\\
  \ge&  \bar J_i(\hat z_{0,T}, a_{0,T}^i) - \epsilon_{1,N}
  \ge  \bar J_i(\hat z_{0,T}, \hat a_{0,T}^i)-\epsilon_{1,N}.
\end{align*}
On the other hand, denoting
$\epsilon_{2,N}=| J_i(\hat a_{0,T}^i, \hat a_{0,T}^{-i})-\bar J_i(\hat z_{0,T}, \hat a_{0,T}^i)|$, we have
 $\lim_{ N\to \infty} \epsilon_{2,N}=0$. Therefore,
$$
J_i(a_{0,T}^i, \hat a_{0,T}^{-i}) \ge J_i(\hat a^i_{0,T}, \hat a_{0,T}^{-i})
-(\epsilon_{1,N} +\epsilon_{2,N}).$$
The theorem follows by taking $\epsilon =\epsilon_{1,N} +\epsilon_{2,N} $.
 \qed

\section{Existence Result}

\label{sec:ex}

Denote  ${\cal Z}_T^{m_0}=\{ z_{0,T}|z_0=m_0, z_t\in [0,1] \mbox{ for } 1\le t\le T\}$.
We introduce the following assumptions.

(H1) $\xi$ has a probability density function denoted by $f_\xi$.

(H2) Consider  the optimal control problem with
cost function
$
\bar J_i(z_{0,T}, a_{0,T}^i)= E\sum_{t=0}^T \rho^t c(x_t^i, z_t, a_t^i).
$
For any $z_{0,T}\in {\cal Z}_T^{m_0}$, there exists $c>0$ such that
the optimal  policy satisfies $a_t^i(x)=a_0$ for all $x\in [0, c]$ and $0\le t\le T$.

We call (H2) the uniformly positive threshold condition for the family of optimal control problems. When the state of the player is small, the effort cost outweighs the extra benefit in further reducing the risk by active control. This holds uniformly with respect to   $z_{0,T}$.

Define the class ${\cal P}_0$ of probability measures  on $\bS$ as follows.
 $\nu\in {\cal P}_0$ if
there exist a constant $c_\nu\ge  0$ and a measurable function $g(x)\ge 0$ defined on $[0,1]$ such that
$$
\nu(B) = \int_B g(x) dx +c_\nu 1_B(0),
$$
where $B\in {\cal B}(\bS)$ and $1_B$ is the indicator function of $B$.
When restricted to $(0, 1]$, $\nu$ is absolutely continuous with respect to the Lebesgue measure $\mu^{\rm Leb}$.

Assume (A1)-(A3) and (H1)-H2) hold, and the distribution of $x_0^i$ is
$\mu_0\in {\cal P}_0$  for this section.

 For given
$z_{0,T}\in {\cal Z}_T^{m_0}$,
let the optimally controlled state process be $x_t^i$ with distribution $\mu_t$.     Define
$w_t= \int_0^1 x \mu_t(dx)$
and the mapping $\Phi$ from $[0,1]^T$ to $[0,1]^T$:
$$
(w_1, \ldots, w_T)= \Phi(z_1, \ldots, z_T).
$$

\begin{lemma} \label{lemma:p}
$\Phi$ is continuous.
\end{lemma}

\proof Let $z_{0,T}\in{\cal Z}_T^{m_0} $ be fixed,  and denote the optimal policy by $a_{0,T}^i$ and the state process by $x_t^i$. Select $ z_{0,T}'\in {\cal Z}_T^{m_0}  $, and denote the corresponding optimal policy by  $ b_{0,T}^{i}$ and the state process by $y_t^{{i}}$.
Let the distribution of $x_t^i$ and $y_t^i$ be  $\mu_t$ and $\mu_t'$, respectively.
Here $\mu_0=\mu_0'\in {\cal P}_0$.
By Lemmas \ref{lemma:y0} and \ref{lemma:y1}, both $\mu_t$ and $\mu_t'$ are in ${\cal P}_0$ for $t\le T$. This ensures that  $\mu_t$ has a small perturbation when the associated positive threshold parameters have a small perturbation.
 By Lemmas \ref{lemma:zza} and \ref{lemma:zzb}, we can first show that
$$
\lim_{z_{0,T}'\to z_{0,T}}\sup_{B\in {\cal B}(\bS)}|\mu_1(B)-\mu_1'(B)|=0.
$$
Repeating the estimate, we further obtain
$$
\lim_{z_{0,T}'\to z_{0,T}}\sup_{B\in {\cal B}(\bS)}|\mu_t(B)-\mu_t'(B)|=0, \quad 0\le t\le T.
$$
Subsequently,
$$
\lim_{z'_{0,T}\to z_{0,T}}\int_0^1 x\mu_t'(dx) = \int_0^1 x\mu_t(dx), \quad 0\le t\le T.
$$
This proves continuity. \qed

\begin{theorem}
There exists a solution $(\hat a^i_{0,T}, \hat z_{0,T})$ to \eqref{dpz}.
\end{theorem}

\proof The theorem follows from Lemma \ref{lemma:p} and  Brouwer's fixed point theorem.
\qed

\section{The Stationary Equations}

\label{sec:se}

\subsection{The stationary form}

Assume (A1)-(A3).
This section introduces a stationary version of \eqref{dpz}.
Take  $z\in \bS$. The value function is independent of time $t$ and so  denoted as $V(x)$.
The dynamic programming equation becomes
$$
V(x)=\min_{a^i} [ c(x, z, a^i)+\rho E[V(x_{t+1}^i)|x_t^i=x]],
$$
which gives
\begin{align}\label{Vsta}
V(x) =  \min \Big[\rho \int_0^1 V( y) Q_0(dy|x) + R(x, z),  \quad \rho V(0) + R(x,z)+ \gamma\Big].
\end{align}
We introduce another equation
\begin{align}
z= \int_0^1x \pi(dx)\label{zxt}
\end{align}
for the probability measure $\pi$.
We say $(\hat z,\hat a^i,\hat \pi)$ is  a stationary solution to \eqref{Vsta}-\eqref{zxt}
if i) the feedback policy $\hat a^i$ is the best response with respect to $\hat z$
in \eqref{Vsta}, ii)  $\{x_t^i, t\ge 0\}$ under the policy $\hat a^i$ has the stationary distribution $\hat \pi$, and
iii) $(\hat z, \hat \pi)$ satisfies \eqref{zxt}.

The equation system \eqref{Vsta}-\eqref{zxt} can be interpreted as follows. For the finite horizon problem, suppose that $T$ is increasing toward $\infty$. If the family of solutions (indexed by different values of $T$) could settle down to a steady-state, we expect for very large $t$, $V(t,x)$ and $z_t$ will be nearly independent of time. This motivates us to introduce
\eqref{Vsta}-\eqref{zxt} as the stationary version of \eqref{dpz}.

\subsection{Value function with general $z$}
\label{sec:sub:val}

Consider a general $z\in \bS$ not necessarily  satisfying \eqref{Vsta}-\eqref{zxt}
simultaneously, and further determine $V(x)$ by \eqref{Vsta}.
Denote $G(x)= \int_0^1 V( y) Q_0(dy|x)$.

\begin{lemma} \label{lemma:LV}
i) Equation \eqref{Vsta} has a unique solution $V\in C([0,1], \mathbb{R})$.

ii) $V$ is strictly increasing.

iii) The optimal policy can be determined as follows:

\quad a) If $\rho G(1)\le \rho V(0)+\gamma$, $a^i(x)\equiv a_0$.

\quad b) If $\rho G(0)\ge \rho V(0)+\gamma$,    $a^i(x) \equiv a_1$.

\quad c) If $\rho G(0)< \rho V(0)+\gamma <\rho G(1)$, there exists a unique $x^*\in (0, 1)$ and $a^i$ is a threshold policy with parameter $x^*$.
\end{lemma}

\proof Part i) will be shown by a fixed point argument.
Define the dynamic programming operator
$$
({\cal L} g)(x) = \min \Big[\rho \int_0^1 g( y) Q_0(dy|x) + R(x, z),
\quad \rho g(0) + R(x,z)+ \gamma\Big],
$$
where $g\in C([0,1], \mathbb{R})$. By the method in proving Lemma
\ref{lemma:con}, it can be shown that ${\cal L}g \in C([0,1], \mathbb{R})$. Now take $g_1, g_2\in  C([0,1], \mathbb{R}).$ Denote
$\hat x= \arg \max |({\cal L}g_2)(x)-({\cal L}g_1)(x)|$. Without  loss of generality, assume $({\cal L}g_2)(\hat x) -({\cal L}g_1)(\hat x) \ge 0$.

Case 1) $({\cal L}g_1)(\hat x)  =\rho \int_0^1 g_1( y) Q_0(dy|\hat x) + R(\hat x, z)$. We obtain
\begin{align*}
0\le ({\cal L}g_2)(\hat x) -({\cal L}g_1)(\hat x) &\le \rho
\int_0^1 (g_2(y)-g_1(y)) Q_0(dy|\hat x)\\
&\le \rho \|g_2-g_1\|.
\end{align*}

Case 2) $({\cal L}g_1)(\hat x) =\rho g_1(0) + R(\hat x,z)+ \gamma$. It follows that \begin{align*}
0\le ({\cal L}g_2)(\hat x) -({\cal L}g_1)(\hat x) &\le \rho |g_2(0)-g_1(0)|\\
& \le \rho \|g_2-g_1\|.
\end{align*}

Combining the two cases, we conclude that ${\cal L}$ is a contraction and has a unique fixed point $V$.

To show ii),
Define $g_{k+1}
= {\cal L} g_k$ for $k\ge 0$, and $g_0=0$.
By the method in Lemma \ref{lemma:in} and  induction, it can be shown that each $g_k$ is increasing on $[0,1]$.
Since $\lim_{k\to \infty}\|g_k-V\|=0$,
$V$ is increasing. Recalling \eqref{Vsta}, we claim that $V$ is strictly increasing.
This proves ii).
 By showing that $G(x)$ is  strictly increasing, we further obtain iii). \qed

For  given $z$,  Lemma \ref{lemma:LV} shows the structure of the optimal policy. Now we specify an optimal policy $a^i(x)$ in terms of  a threshold parameter $\theta(z)$ by the  following rule.
i) If $\rho V(0)+\gamma \le \rho G(0)$, then $a^i(x)\equiv a_1$
with $\theta(z)=0$;
 ii) if $\rho G(0)<\rho V(0)+\gamma <\rho G(1)$,  $a^i(x)$ is a threshold policy with
 $\theta(z)\in (0, 1)$;
iii) if $\rho V(0)+\gamma =\rho G(1)$, then $a^i(x)$  has  $\theta(z)=1$;
iv) if $\rho V(0)+\gamma >\rho G(1)$,
 $a^i(x)\equiv a_0$ for which we  formally denote $\theta(z)=1^+$.

\begin{remark}
 For the case  $\rho V(0)+\gamma =\rho G(1)$, the above rule gives
 $a^i(1)=a_1$ and $a^i(x)=a_0$ for $x<1$, which is slightly different from Lemma \ref{lemma:LV} but still attains optimality.
\end{remark}

%If the resulting optimal policy is $a^i(x)\equiv a_0$, we formally denote %it as a threshold policy with parameter $\theta(z)= 1^+$;
%otherwise, it is a usual threshold policy with parameter $\theta(z)\in [0,1]$.

\subsection{Stationary distribution for a given threshold policy}

Suppose that $a^i$ is a threshold policy with parameter $\theta\in (0, 1)$. Denote the corresponding state process by $\{x_t^{i,\theta}, t\ge 0\}$, which is a Markov process. Let the probability measure $P^t(x,\cdot)$ on ${\cal B}(\bS)$ be the distribution of $x_t^{i,\theta}$ given $x_0^{i,\theta}=x\in \bS$.

We introduce a further condition on $\xi$.

(A4) $\xi$ has a probability density function $f_\xi(x)>0$ a.e. on $\bS$.

\begin{theorem} \label{theorem:erg}
For $\theta\in (0, 1)$,  $\{x_t^{i,\theta}, t\ge 0\}$
is uniformly ergodic  with stationary probability distribution $\pi_\theta$, i.e.,
\begin{align}
\sup_{x\in \bS}\|P^t(x, \cdot)
-\pi_\theta\|_{\rm TV}\le K r^t \label{ppit}
\end{align}
for some constants $K>0$ and $r\in (0, 1)$,
where $\|\cdot\|_{\rm TV}$ is the total variation norm of signed measures.
 %the probability measure $P^t(x,\cdot)$  converges to a stationary distribution  %$\pi_\theta$ in the total variation norm, uniformly with respect to  $x$.
\end{theorem}

\proof See appendix B.    \qed

%\newpage

%\subsection{The fixed point problem}

%For \eqref{Vsta}-\eqref{zxt},
%we  look for a solution $(z, a^i,\pi)$ from the class ${\cal C}$ of solutions  where $z\in\bS$ and $a^i$ is a %threshold policy with parameter $\theta \in [0,1]$ or $\theta=1^+$.
%The existence analysis of  \eqref{Vsta}-\eqref{zxt}  reduces to a fixed point problem.
%We start with any $z\in \bS$. If $\theta(z)\in (0, 1)$, let the stationary distribution
%of $x_t^{i, \theta(z)}$ be $\pi_{\theta(z)}$ and define
%$$
%\Psi(z)= \int_0^1 x \pi_{\theta(z)}(dx).
%$$
%If $\theta(z)=0$, we have $E x_t^{i, \theta(z)}=0$ for all $t\ge 1$, and define $\Psi(z)=0$.
%If $\theta(z)=1$ or $\theta(z)=1^+$, define $\Psi(z)=1$.
%The fixed point equation becomes
%$z=\Psi(z).$

%By considering a product form $R(x,z)= R_1(x)R_2(z)$ and assuming that both
%$R_1$ and $R_2$ are increasing, we can show that  \eqref{Vsta}-\eqref{zxt}
%has at most one solution in ${\cal C}$. An increasing function $R_2$ means positive externalities of
%the mean field. For detailed analysis, see \cite{HM16}.

\subsection{Comparison theorems}

Denote
$z(\theta) = \int_0^1 x \pi_\theta(dx).$
We have the first comparison  theorem on  monotonicity.
\begin{theorem}\label{theorem:zmono}
 $z(\theta_1)\le z(\theta_2)$ for $0<\theta_1<\theta_2<1$.
\end{theorem}

\proof See appendix D. \qed

In the further analysis,
we consider the case where $R$ takes the product form
$
R(x, z)= R_1( x) R_2(z),
$
and where $R$ still satisfies (A2) and $R_1\ge 0$, $R_2> 0$. We further assume

(A5) $R_2>0$ is strictly increasing on $\bS$.

This assumption indicates positive externalities since an individual  benefits from  the decrease of the population average state. This condition has a crucial role in the uniqueness analysis.

Given the product form of $R$, now \eqref{Vsta} takes the form
\begin{align*}
V(x) =  &\min \Big[\rho \int_0^1 V( y) Q_0(dy|x) + R_1(x)R_2(  z),
\quad \rho V(0) + R_1(x)R_2( z)+ \gamma\Big].
\end{align*}
%If the optimal policy is $a^i(x)\equiv 0$, we formally denote it as a threshold policy with parameter %$\theta(z)= 1^+$.
%Otherwise, it is a usual threshold policy with parameter $\theta(z)\in [0,1]$.
Consider $0\le z_2<z_1\le 1$ and
\begin{align}
V_l(x) =  &\min \Big[\rho \int_0^1 V_l( y) Q_0(dy|x) + R_1(x)R_2(  z_l),  \quad \rho V_l(0) + R_1(x)R_2( z_l)+ \gamma\Big]. \label{Vlr12}
\end{align}
Denote the optimal policy as a threshold policy  with parameter $\theta_l$ in $[0, 1]$ or equal to $1^+$, where we follow the rule in Section \ref{sec:sub:val}
to interpret $\theta_l=1^+$.
We state the second comparison theorem about  the threshold parameters under different mean field parameters $z_l$.

\begin{theorem} \label{theorem:s1s2}
 $\theta_1$ and $\theta_2$ in \eqref{Vlr12} are
specified according to the following scenarios:

i) If $\theta_1=0$, then we have either $\theta_2\in [0, 1]$ or
$\theta_2=1^+$.

ii)   If $\theta_1\in (0,1)$, we have either  a)
$\theta_2\in (\theta_1, 1)$, or b) $\theta_2=1 $,
or c) $\theta_2= 1^+$.

iii) If $\theta_1=1$, $\theta_2=1^+$.

iv) If $\theta_1=1^+$, $\theta_2= 1^+$.
\end{theorem}

\proof Since $R_2(z_1)> R_2(z_2)>0$, we divide both sides of \eqref{Vlr12} by $R_2(z_l)$ and define
$
\gamma_l= \frac{\gamma}{R_2(z_l)}.
$
Then $0<\gamma_1<\gamma_2$. The dynamic programming equation reduces to \eqref{vlr12}.
Subsequently, the optimal policy is determined according to Lemma \ref{lemma:interv}. \qed

\subsection{Uniqueness}

We look for a solution $(z, a^i,\pi)$ from the class ${\cal C}$ of solutions  where $z\in\bS$ and $a^i$ is a threshold policy with parameter $\theta \in [0,1]$ or $\theta=1^+$.

\begin{theorem}
Under (A1)-(A5) with $R(x,z)=R_1(x)R_2(z)$, the equation system \eqref{Vsta}-\eqref{zxt} has at most one solution in  ${\cal C}$.
\end{theorem}

\proof
Assume two different solutions
\begin{align}
(z_1, a^i, \pi)\neq (z_2, b^i, \nu).\label{zab}
\end{align}
If $z_1=z_2$, \eqref{Vsta} ensures $a^i=b^i$, and subsequently $\pi=\nu$. This is a contradiction to two different solutions.
Now we can assume
\begin{align}0\le z_2<z_1\le 1. \label{z211}
\end{align}
We check different scenarios listed in Theorem \ref{theorem:s1s2}.
If  $\theta_1\in (0, 1)$ so that $\theta_2\in (\theta_1, 1)$, Theorem
\ref{theorem:zmono} implies $z_1\le z_2$, which contradicts \eqref{z211}.
For all remaining scenarios, it is easy to show $z_1\le z_2$, which again contradicts \eqref{z211}.
Therefore,  assumption  \eqref{zab} does not hold. Uniqueness follows. \qed

\section{Conclusion}

This paper considers mean field games in a framework of multi-agent Markov decision processes (MDP).
Each player has a monotone cost function and can apply resetting control to the state process.
Decentralized strategies are obtained as threshold policies. We further examine a system of stationary equations of the mean field game and study uniqueness of the solution under positive externalities.

\label{sec:con}

\section*{Appendix A:  Technical Lemmas for Section \ref{sec:ex}}
\renewcommand{\theequation}{A.\arabic{equation}}
\setcounter{equation}{0}
\renewcommand{\thetheorem}{A.\arabic{theorem}}
\setcounter{theorem}{0}
\renewcommand{\thelemma}{A.\arabic{lemma}}
\setcounter{lemma}{0}

Let $X$ be a random variable with distribution $\nu\in {\cal P}_0$.
Set $x_t^i=X
$. Define $Y_0=x_{t+1}^i$ by applying $a_t^i\equiv a_0$.
Further define $Y_1=x_{t+1}^i$ by applying the threshold policy $a_t^i$ with parameter $r\in (0, 1)$.
 Then
\begin{align}
P(Y_0\in B) =\int_0^1 Q_0(B|x) \nu(dx), \quad B\in {\cal B}(\bS).
\end{align}

\begin{lemma}\label{lemma:y0}
The distribution of  $Y_0$ is in ${\cal P}_0$.
\end{lemma}

\proof We can directly show that the probability density function of $Y_0$ is
$$
g(y) = \int_{0\le x<y} \frac{1}{1-x}  f_\xi \left(\frac{y-x}{1-x}\right) \nu (dx), \quad y\in (0,1).
$$
 In this case $P(Y_0=0)=0$. \qed

\begin{lemma}\label{lemma:y1}
The distribution of $Y_1$ is in ${\cal P}_0$.
\end{lemma}

\proof It is clear that
$ P(Y_1=0)= P(X\ge r).$
 The distribution of $Y_1$ restricted on $(0, 1]$ is absolutely continuous with respect to $\mu^{\rm Leb}$. Denote
$$
g(y) = \int_{0\le x<y\wedge r} \frac{1}{1-x}  f_\xi \left(\frac{y-x}{1-x}\right) \nu (dx).
$$
 Then for $B\in {\cal B}(\bS)$,
$$
P(Y_1\in B)= \int_B g(y) dy + P(X\ge r) 1_B(0).
$$
The lemma follows. \qed

Let $z_{0,T}\in {\cal Z}_T^{m_0}$ be fixed, and denote the associated optimal policy by $ a_{0,T}^i$.
Select another sequence $z_{0,T}'\in {\cal Z}_T^{m_0}  $  and let the optimal policy be denoted by $ b_{0,T}^i$. Denote $d(z_{0,T}',z_{0,T})= \sum_{k=1}^T|z_k'-z_k|.$ Write $z_{0,T}'\to z_{0,T}$ when $d(z_{0,T}',z_{0,T})\to 0.$ Fix any $t\le T-1$. Based on (H2),  we consider two cases.

Case A) $a_t^i(x) =a_0$  for all $x\in \bS$.

\begin{lemma} \label{lemma:zza}
For case A) and any $\epsilon>0$,
there exists $\delta>0$ such that for all $z_{0,T}'$ satisfying
  $d(z_{0,T}',z_{0,T})\le \delta$, we have either

i) $ b_t^i(x) = a_0 $ for all $x\in \bS$, or

 ii) there exists $r'\in (0,1)$ such that
$ b_t^i $ is a threshold policy with parameter $r'$ and $0<1-r'<\epsilon $.
\end{lemma}

\proof
When $ z_{0,T}' \in {\cal Z}_T^{m_0} $ is used in the optimal control problem in (H2),  denote the value function by
$\acute{V}(t,x)$. Define $\acute{G}_t(x)$ in place of $G_t(x)$. Then $\acute{G}_t(x)$ is continuous and strictly increasing for each $t$.
Since $V$ depends on $z_{0,T}$ continuously,
\begin{align}
\lim_{z'_{0,T}\to z_{0,T}} \sup_{x\in \bS}|\acute{V}(t,x)-V(t,x)|=0.
\label{vcon}
\end{align}
 Consider any $0<\epsilon<1$.
We only need to treat the following two scenarios.

1) $\rho G_{t+1}(1) < \rho V(t+1, 0)+\gamma$.

By \eqref{vcon}, there  exists $\delta>0$ such that for all $z_{0,T}'$ satisfying
  $d(z_{0,T}',z_{0,T})\le \delta$, we have
\begin{align}
\rho \acute{G}_{t+1}(1) < \rho \acute{V}(t+1, 0)+\gamma. \label{gvp}
\end{align}
Then $ b_t^i(x)=a_0 $ for all $x\in \bS$. So i) holds.

2) $\rho G_{t+1}(1) = \rho V(t+1, 0)+\gamma.$

 Then
\begin{align}
\rho G_{t+1}(1-\epsilon) < \rho V(t+1, 0)+\gamma. \label{g1e}
\end{align}
By \eqref{g1e}, there exists $\delta>0$ such that for all $z_{0,T}'$ satisfying
  $d(z_{0,T}',z_{0,T})\le \delta$, we have
\begin{align}
\rho \acute{G}_{t+1}(1-\epsilon) < \rho \acute{V}(t+1, 0)+\gamma. \label{gpvpg}
\end{align}
For such  $z_{0,T}'$, if
 $\rho \acute{G}_{t+1}(1) \le \rho \acute{V}(t+1, 0)+\gamma,$
    we select $ b^i_t(x)=a_0$ for all $x\in\bS$ and then i) holds.
 If $z_{0,T}'$ results in
\begin{align}
\rho \acute{G}_{t+1}(1) > \rho \acute{V}(t+1, 0)+\gamma, \nonumber %\label{gpga}
\end{align}
 by \eqref{gpvpg}, we can find $r'\in (1-\epsilon,1)$ such that
 \begin{align}
\rho \acute{G}_{t+1}(r') = \rho \acute{V}(t+1, 0)+\gamma, \nonumber  %\label{gprp}
\end{align}
which further determines $b_t^i$ as a threshold policy with parameter $r'$.
 \qed

Case B). There exists $r\in (0, 1)$ such that $ a_t^i$
is a threshold policy with parameter $r$.

\begin{lemma}\label{lemma:zzb}
For case B), when $d(z'_{0,T},z_{0,T})$ is sufficiently small, $b_t^i$ is a threshold policy with parameter
$r'\in (0, 1)$ and in addition,
$r'\to r$ as $z'_{0,T}\to z_{0,T}$.
\end{lemma}

\proof We have
$
\rho G_{t+1}(r)= \rho V(t+1,0)+\gamma.
$
Fix  a small $\epsilon>0$ such that $(r-\epsilon, r+\epsilon)\subset (0 ,1)$. Then
$$
\rho G_{t+1}(r-\epsilon)<\rho V(t+1, 0)+\gamma, \quad
\rho G_{t+1}(r+\epsilon)>\rho V(t+1, 0)+\gamma.
$$
By \eqref{vcon}, we may select  $\delta>0$  such that for all $z_{0,T}'$ satisfying $d(z'_{0,T}, z_{0,T})\le \delta$, we have
$$
\rho \acute{ G}_{t+1}(r-\epsilon)<\rho \acute{V}(t+1, 0)+\gamma, \quad
\rho \acute{G}_{t+1}(r+\epsilon)>\rho \acute{V}(t+1, 0)+\gamma.
$$
Since $\acute{G}_{t+1}(x)$ is strictly increasing, there exists a unique
$r'\in (r-\epsilon, r+\epsilon)$ such that
$$
\rho \acute{G}_{t+1}(r')=\rho \acute{V}(t+1, 0)+\gamma,
$$
where $r'$ depends on $ z_{0,T}'$.
This in turn determines the threshold policy $b_t^i$ with parameter $r'$. Since $\epsilon$ can be arbitrarily small, the last part of the lemma follows. \qed

\section*{Appendix B. Proof of Theorem \ref{theorem:erg}}
\renewcommand{\theequation}{B.\arabic{equation}}
\setcounter{equation}{0}
\renewcommand{\thetheorem}{B.\arabic{theorem}}
\setcounter{theorem}{0}
\renewcommand{\thelemma}{B.\arabic{lemma}}
\setcounter{lemma}{0}

Consider $0<\theta<1$.
The definitions of irreducibility, aperiodicity and a small set
follow those in
\cite{MT09}.
Let $\delta_x$ be the dirac measure at  $x\in \mathbb{R}$.
Let $\varphi: =\delta_0$. So $\delta_0(B)=1_B(0)$ for $B\in {\cal B}(\bS)$.

Throughout this appendix, we write $x_t:= x_t^{i,\theta}$ in order to keep the notation light.

\begin{lemma}
$\{x_t, t\ge 0\}$ is $\varphi$-irreducible.
\end{lemma}

\proof We can directly verify that
\begin{align*}
&P(x_2=0|x_0=x)>0, \quad x\in [0,\theta), \\
& P(x_1=0|x_0=x)=1, \quad x\in [\theta, 1].
\end{align*}
The above probabilities of the process are calculated by setting the  distribution for $x_0$ as the dirac measure ${\bf \delta}_x$. This implies that $\{x_t, t\ge 0\}$ is $\varphi$-irreducible.
 \qed

\begin{lemma}
$\{x_t, t\ge 0\}$ is aperiodic.
\end{lemma}

\proof Define $C_s= \{0\}$. Denote $\epsilon_0=\int_\theta^1 f_\xi (y) dy >0$ and the measure $\nu=\epsilon_0 \delta_0$. Then
\begin{align*}
P(x_2=0|x_0=0 )& \ge P(x_2=0, x_1\ge \theta|x_0=0)\\
   & = P(x_1\ge \theta |x_0=0) \\
   %& =P(\xi \ge \theta)\\
   & =\epsilon_0.
\end{align*}
For any $B\in {\cal B}(\bS)$, then
\begin{align}
P(x_2\in B|x_0=0)\ge \nu (B). \label{sm2}
\end{align}
Therefore, we can take $C_s$ as a small set with $\nu(C_s)=\epsilon_0$.
Given $x_0=0\in C_s$, we further check
\begin{align*}
P(x_3=0|x_0=0) &\ge P(x_3=0, x_2\ge \theta, x_1<\theta|x_0=0)\\
    &= P(x_2\ge \theta, x_1<\theta|x_0=0).
\end{align*}
Let $ \xi,\xi_1, \xi_2$ be i.i.d. random variables. Then
\begin{align*}
 P(x_2\ge \theta, x_1<\theta|x_0=0)
 &=P(\xi_1+(1-\xi_1)\xi_2\ge \theta, \xi_1<\theta)\\
 &\ge P(\xi_2\ge \theta, \xi_1<\theta).
\end{align*}
Hence
\begin{align*}
P(x_3=0|x_0=0)\ge
\int_0^\theta f_\xi (y) dy \int_\theta^1 f_\xi (y)dy.
\end{align*}
Denote $ \epsilon_1=\int_0^\theta f_\xi (y) dy $.
Then for any $B\in {\cal B}(\bS)$,
 \begin{align}
P(x_3\in B|x_0=0)\ge \epsilon_1\nu(B). \label{sm3}
\end{align}
Since time indices 2 and 3 in \eqref{sm2} and \eqref{sm3} have the greatest common divisor equal to 1, $\{x_t, t\ge 0\}$ is aperiodic \cite[pp. 112-114]{MT09}. \qed

The Markov process $\{x_t, t\ge 0\}$ is said to satisfy Doeblin's condition if there exist a probability measure $\phi$ on ${\cal B}(\bS)$ and $\epsilon <1$, $\eta>0$, $m\ge 0$, such that $\phi(B)>\epsilon$ implies
$$
\inf_{x\in \bS} P(x_m\in B |x_0=x)\ge \eta.
$$

\begin{lemma}
Doeblin's condition holds for $\{x_t, t\ge 0\}$.
\end{lemma}

\proof
We take $\phi=\delta_0$, and in this case $\phi(B)>0$ implies $ 0\in B $. It suffices to show
$$
\inf_{x\in \bS} P(x_4=0 |x_0=x)\ge \eta .
$$
For $x\in [\theta, 1]$,
\begin{align}
P(x_4=0|x_0=x)& \ge P(x_4=0, x_3\ge \theta, x_2<\theta, x_1=0|x_0=x)\nonumber \\
  % &=P(x_4=0, x_3\ge \theta, x_2<\theta, x_1=0|x_0=x) \nonumber \\
   &= P(x_3=0, x_2\ge \theta, x_1<\theta|x_0=0) \nonumber\\
   &= P(x_2\ge \theta, x_1<\theta|x_0=0)\nonumber  \\
   &\ge \epsilon_0\epsilon_1. \label{db1}
\end{align}
Let $\xi$, $\xi_k$, $k=1,2,3$, be i.i.d.
For $x\in [0, \theta)$,
\begin{align}
P(x_4=0|x_0=x) &\ge P(x_4=0, x_3\ge \theta, x_2=0, x_1\ge \theta|x_0=x) \nonumber \\
   & = P( x_3\ge \theta, x_2=0, x_1\ge \theta|x_0=x) \nonumber \\
   &= P(\xi_3\ge \theta) P(x_2=0, x_1\ge \theta|x_0=x)\nonumber  \\
   &= P(\xi_3\ge \theta) P(x+(1-x) \xi_1\ge \theta)\nonumber  \\
   &\ge P(\xi_3\ge \theta) P(\xi_1\ge \theta)=\epsilon_0^2 .  \label{db2}
\end{align}
By \eqref{db1} and \eqref{db2}, Doeblin's condition holds with $\epsilon=\frac12$, $m=4$, $\eta=\epsilon_0^2\epsilon_1>0$. \qed

{\it Proof of Theorem \ref{theorem:erg}.}
Since $\{x_t, t\ge 0\}$ is aperiodic and satisfies
Doeblin's condition, by \cite[pp. 394, Theorem 16.0.2]{MT09}, \eqref{ppit} holds.
\qed

\section*{Appendix C: An Auxiliary MDP}
\renewcommand{\theequation}{C.\arabic{equation}}
\setcounter{equation}{0}
\renewcommand{\thetheorem}{C.\arabic{theorem}}
\setcounter{theorem}{0}
\renewcommand{\thelemma}{C.\arabic{lemma}}
\setcounter{lemma}{0}
\renewcommand{\theremark}{C.\arabic{remark}}
\setcounter{remark}{0}

%\subsection{}

Assume (A1)-(A4).
This appendix introduces an auxiliary optimal control problem
to show the effect of the effort cost on the threshold parameter of the optimal policy.
The state and control processes
$\{(x_t^i, a_t^i), t\ge 0 \}$  are specified by  \eqref{xa0}-\eqref{xa1}. The cost has the form
\begin{align}
J_i^r=E\sum_{t=0}^\infty \rho^t (R_1(x^i_t)+  r 1_{\{a_t^i=a_1\}} ),
\label{jir}
\end{align}
where $R_1$ is continuous and strictly increasing on $[0, 1]$ and $r\in (0, \infty)$.
Let $r$ take two different values  $0<\gamma_1<\gamma_2$ and write the
corresponding dynamic programming equation
\begin{align}
v_l(x) =\min \left\{\rho \int_0^1 v_l (y) Q_0(dy|x)+R_1(x), \quad \rho v_l(0) +R_1(x)+\gamma_l   \right\} , \quad l=1,2,\ x\in \bS.\label{vlr12}
\end{align}

By the method in proving Lemma \ref{lemma:LV}, it can be shown that there exists a unique solution
$v_l\in C([0,1], \mathbb{R})$ and that the optimal policy $a^{i,l}(x)$ is a threshold policy.
If
$ \rho \int_0^1 v_l (y) Q_0(dy|1)<\rho v_l(0)+\gamma_l,
$
$a^{i,l}(x)\equiv  a_0$, and we follow the notation in Section \ref{sec:sub:val}
to formally denote it as a threshold policy with parameter $\theta_l=1^+$. Otherwise, $a^{i,l}(x)$ is a
 threshold policy with parameter $\theta_l\in [0,1]$, i.e., $a^{i,l}(x)= a_1$ if $x\ge \theta_l$, and $a^{i,l}(x)=a_0$ if $x<\theta_l$.

\begin{lemma}\label{lemma:the12}
 If $\theta_1\in (0, 1)$, $\theta_2\ne \theta_1$.
\end{lemma}

\proof
We prove by contradiction. Suppose for some $\theta\in (0,1)$,
\begin{align}
\theta_1=\theta_2=\theta.  \label{th12}
\end{align}
Under  assumption \eqref{th12},
 the resulting optimal policy leads to  the representation (see e.g. \cite[pp. 22]{H93})
$$
v_l(x)= E\sum_{t=0}^\infty \rho^t \left[R_1( x_t^{i})+\gamma_l 1_{\{ a_t^i=a_1\}}\right], \quad l=1,2,
$$
where $\{ x_t^i, t\ge 0\}$ is generated by the threshold policy $a_t^i( x_t^i)$ with parameter $\theta$ and $x_0^i=x$. Denote $\delta_{21}=\gamma_2-\gamma_1$.

For any fixed $x\ge \theta$ and $ x_0^i=x$, denote the resulting optimal state and control processes by $\{(\hat x_t^i, \hat a_t^i), t\ge 0  \}$.
Then $\hat a_0^i =a_1$ w.p.1.,  and
$$
v_2(x)-v_1(x)= \delta_{21} +\delta_{21} E\sum_{t=1}^\infty \rho^t1_{\{ \hat a_t^i=a_{1}\}}, \qquad x\ge \theta.
 $$

Next consider  $x_0^i=0$ and denote the optimal state and control processes by
$\{(\check{x}_t^i, \check{a}_t^i), t\ge 0  \}$.
Then
$$
v_2(0)-v_1(0)= \delta_{21} E\sum_{t=0}^\infty \rho^t1_{\{ \check{a}_t^i=a_{1}\}}
=: \Delta.
$$

It is clear that $\hat x_1^i=0$ w.p.1. By the optimality principle,
$\{(\hat x_t^i, \hat a_t^i), t\ge 1  \}$ may be interpreted as the optimal state and control processes of the MDP with initial
state 0 at $t=1$. Hence the two processes $\{(\hat x_t^i, \hat a_t^i), t\ge 1  \}$ and $\{(\check{x}_t^i, \check{a}_t^i), t\ge 0  \}$, where ${\check x}_0^i=0$, have the same finite dimensional distributions. In particular, $\hat a_{t+1}^i$ and $\check{a}_t^i$ have the same distribution for $t\ge 0$.
Therefore,
$$
E\sum_{t=1}^\infty \rho^{t-1}1_{\{ \hat a_t^i=a_{1}\}}
=
E\sum_{t=0}^\infty \rho^t1_{\{ \check{a}_t^i=a_{1}\}}.
$$
It follows that
\begin{align}
v_2(x)-v_1(x)= \delta_{21}+\rho \Delta, \qquad \forall x\ge  \theta. \label{v12d}
\end{align}
Combining
  \eqref{vlr12} and \eqref{th12} gives
 \begin{align*}
\rho \int_0^1 v_l(y) Q_0(dy|\theta)=  \rho v_l(0)+\gamma_l  ,
 \qquad l=1, 2,
 \end{align*}
 which implies
 \begin{align}
\rho \int_0^1 [v_2(x)-v_1(x)] Q_0(dx|\theta)  =\delta_{21}+ \rho \Delta.  \label{v21del}
 \end{align}
By  $Q_0([0, \theta)|\theta)=0$ and \eqref{v12d},  \eqref{v21del} further yields
 \begin{align*}
 \rho( \delta_{21}+ \rho \Delta)
=  \delta_{21}+\rho \Delta,
 \end{align*}
 which is impossible since $0<\rho<1$ and $\delta_{21}+\rho \Delta>0$.
 Therefore, \eqref{th12} does not hold. This completes the proof.
  \qed

%\subsection{ $r$ on the Threshold}

For the MDP with cost
\eqref{jir},
we continue to analyze the  dynamic programming equation
\begin{align}\label{Vstab}
v_r(x) =  \min \Big[\rho \int_0^1 v_r( y) Q_0(dy|x) + R_1(x),  \quad \rho v_r(0) + R_1(x)+ r\Big].
\end{align}
For each fixed $r \in (0, \infty )$, we obtain the optimal policy as a threshold policy with parameter $\theta(r)$.
 By evaluating the cost \eqref{jir}   associated with  the two policies
 $a_t^i(x_t^i)\equiv a_0$ and $a_t^i(x_t^i)\equiv a_1$, respectively,  we have the prior estimate
\begin{align}
 v_r(x)\le \min\left\{ \frac{R_1(1)}{1-\rho },\ R_1(x)+\frac{r +\rho R_1(0) }{1-\rho}   \right\}. \label{rvr}
\end{align}
On the other hand, let $\{x_t^i, t\ge 0\}$ with $x_0^i=x$ be generated by any fixed Markov policy. Then
\begin{align}
E\sum_{t=0}^\infty \rho^t (R_1(x_t^i)+ r 1_{\{a_t^i=a_1\}})\ge R_1(x) +
\sum_{t=1}^\infty \rho^t R_1(0), \nonumber
\end{align}
which implies
\begin{align}
v_r(x)\ge R_1(x)+\frac{\rho R_1(0)  }{1-\rho }. \label{vrge}
\end{align}

If $r>\frac{\rho R_1(1)}{1-\rho }$, it follows from \eqref{rvr} that
\begin{align}
\rho \int_0^1 v_r( y) Q_0(dy|x) < \rho v_r(0)+r,   \quad \forall x,  \label{th1}
\end{align}
i.e., $\theta(r)=1^+$.

\begin{lemma}\label{lemma:th0}
There exists $\delta>0$ such that for all $0<r<\delta$,
\begin{align}
\rho \int_0^1 v_r( y) Q_0(dy|x)> \rho v_r(0)+r, \quad \forall x,  \label{th0}
\end{align}
and so $\theta(r)=0$.
\end{lemma}

\proof By \eqref{vrge},
\begin{align}
\rho \int_0^1 v_r( y) Q_0(dy|x) &\ge \rho\int_0^1 R_1(y) Q_0(dy|x)
+\frac{\rho^2 R_1(0)}{1-\rho } \nonumber \\
 & \ge \rho\int_0^1 R_1(y) Q_0(dy|0)
+\frac{\rho^2 R_1(0)}{1-\rho }  , \nonumber
\end{align}
and \eqref{rvr} gives
$$
\rho v_r(0)+r \le
 %\rho  R_1(0) +\frac{\rho ( r +\rho R_1(0)) }{1-\rho }+r =
\frac{ \rho R_1(0)}{1-\rho }+
\frac{ r }{1-\rho }.
$$
Since $R_1(x)$ is strictly increasing,
$$
C_{R_1}: =\int_0^1 R_1(y) Q_0(dy|0)- R_1(0)>0.
$$
and
$$
\rho \int_0^1 v_r( y) Q_0(dy|x)- (\rho v_r(0)+r) \ge \rho C_{R_1}- \frac{r}{1-\rho}.
$$
It suffices to take $\delta =\rho (1-\rho) C_{R_1}.$ \qed

Define the nonempty sets
$$
{\cal R}_{a_0}=\{r>0| \eqref{th1} \mbox{ hods}\}, \quad {\cal R}_{a_1}=\{r>0|  \eqref{th0} \mbox{ holds}\}  .
$$

\begin{remark}\label{reC:int}
We have $(\frac{\rho R_1(1)}{1-\rho }, \infty   )\subset {\cal R}_{a_0} $ and  $(0, \delta)\subset {\cal R}_{a_1} $.
\end{remark}

\begin{lemma} \label{lemma:connect}
 Let $(r, v_r)$ be the parameter and the associated solution in \eqref{Vstab}.

i) If $r>0$ satisfies
\begin{align}
\rho \int_0^1 v_r(y) Q_0(dy|x) \le \rho v_r(0) +r, \quad \forall x,
\label{rth1}
\end{align}
  then  any $r'>r$ is in ${\cal R}_{a_0}$.

ii) If $r>0$ satisfies
\begin{align}
\rho \int_0^1 v_r(y) Q_0(dy|x) \ge \rho v_r(0) +r,\quad \forall x,
\label{rth0}
\end{align}
then any  $r'\in (0,r)$ is in ${\cal R}_{a_1}$.
\end{lemma}

\proof i) For $r'>r$, $v_{r'}$ is uniquely solved from \eqref{Vstab} with $r'$ in place of $r$. We can use \eqref{rth1} to  verify
$$
v_r(x)= \min \left[\rho \int_0^1 v_r(y) Q_0(dy|x) +R_1(x) , \quad \rho v_r(0) +R_1(x)+r'\right].
$$
Hence $v_{r'}=v_r$ for all $x\in [0,1]$. It follows that
$\rho \int_0^1 v_{r'}(y) Q_0(dy|x) <  \rho v_{r'}(0)+r'$ for all $x$. Hence  $r'\in {\cal R}_{a_0}$.

ii) By \eqref{Vstab} and  \eqref{rth0},
\begin{align*}
 v_r(0) = \frac{ R_1(0)+r}{1-\rho }  , %\label{vr0r}
 \end{align*}
 and subsequently,
\begin{align*}
 v_r(x) &= \rho v_r(0) +R_1(x) +r \nonumber  \\
     &= \frac{\rho R_1(0)+r}{1-\rho }  +R_1(x).  %\label{vrxsol}
\end{align*}
By substituting $v_r(0)$ and $v_r(x)$ into \eqref{rth0}, we obtain
\begin{align}
\rho R_1(0) +r \le \rho \int_0^1 R_1(y) Q_0(dy|x), \quad \forall x.
\label{rhorrle}
\end{align}
Now for $0<r'<r$, we construct $v_{r'}(x)$, as a candidate solution to
\eqref{Vstab} with $r$ replaced by $r'$,    to satisfy
\begin{align}
v_{r'}(0)= \rho v_{r'}(0)+R_1(0)+r', \quad
v_{r'}(x) = \rho v_{r'}(0) +R_1(x) +r', \label{defvp}
\end{align}
which gives
\begin{align}
v_{r'}(x) =\frac{\rho R_1(0)+r'  }{1-\rho }  +R_1(x).
\label{ss}
\end{align}

We  show that $v_{r'}(x)$ in \eqref{ss} satisfies
\begin{align}
\rho v_{r'}(0) +r' <\rho \int_0^1 v_{r'}(y) Q_0(dy|x), \quad \forall x,
\label{rovrp}
\end{align}
which is equivalent to
$$
\rho R_1(0) +r' <\rho \int_0^1 R_1(y) Q_0(dy|x), \quad \forall x,
 $$
which in turn follows from \eqref{rhorrle}.
By \eqref{defvp} and  \eqref{rovrp},  $v_{r'}$ indeed satisfies
\eqref{Vstab} with $r$ replaced by $r'$.  So $r'\in {\cal R}_{a_1}$. \qed

%Still let
%$(v_r, r)$ be the solution and parameter pair in \eqref{Vstab}. Denote the %sets
%$$
%\overline \Gamma_0=\{r> 0| \eqref{rth0} \mbox{ holds} \}, \quad \overline %\Gamma_1=\{r> 0| \eqref{rth1} \mbox{ holds} \}.
%$$

Further define
$$
\underline{r} =\sup {\cal R}_{a_1}, \quad \overline r= \inf  {\cal R}_{a_0}.
$$

\begin{lemma} \label{lemma:interv}
i) $\underline r$ satisfies
\begin{align}
\rho \int_0^1 v_{\underline r}(y) Q_0(dy|0) = \rho v_{\underline r}(0) +
\underline {r}, \nonumber
\end{align}
and $\theta(\underline r)=0$.

ii)
$\overline r$ satisfies
\begin{align}
\rho \int_0^1 v_{\overline r}(y) Q_0(dy|1) =\rho v_{\overline r}(1) = \rho v_{\overline r}(0) +
\overline {r}. \nonumber
\end{align}
and $\theta(\overline r)=1$.

iii)
We have  $0<\underline r<\overline r<\infty$.

iv) $\theta(r)$  is continuous and strictly increasing on $[\underline r, \overline r]$.
\end{lemma}

\proof i)-ii) By Lemmas \ref{lemma:th0} and
\ref{lemma:connect}, we have
$ 0<\underline r\le \infty$ and $0\le  \overline r<\infty$.
Assume $\underline r=\infty$; then ${\cal R}_{a_1} =(0,\infty)$ giving
 ${\cal R}_{a_0}=\emptyset$, a contradiction. So $0<\underline r<\infty$.
For $\delta >0$ in Lemma \ref{lemma:th0}, we have $(0, \delta)\subset
{\cal R}_{a_1}$. Therefore, $0<\bar r <\infty$.
Note that $v_r$ depends on the parameter $r$ continuously, i.e.,
 $\lim_{|r'-r|\to 0}\sup_x|v_{r'}(x)-v_r(x)|=0$.
Hence
$$
\rho \int_0^1 v_{\underline r}(y) Q_0(dy|0) \ge  \rho
v_{\underline r}(0) +
\underline {r}.
$$
Now assume
\begin{align}
\rho \int_0^1 v_{\underline r}(y) Q_0(dy|0) >  \rho v_{\underline r}(0) +
\underline {r}. \label{rhor16}
\end{align}
Then there exists a sufficiently small $\epsilon>0$ such that \eqref{rhor16} still holds when
$(\underline r+\epsilon , v_{\underline r+\epsilon}   )$ replaces
$({\underline r} , v_{\underline r} )$; since $g(x)=\int_0^1
 v_{\underline r+\epsilon}(y) Q_0(dy|x) $ is increasing in $x$,  then $\underline r+\epsilon \in {\cal R}_{a_1}$, which is impossible. Hence \eqref{rhor16} does not hold, and this proves i). ii) can be shown in a similar manner.

To show iii),
  %It is clear that ${\cal R}_{a_0}\cap {\cal  R}_{a_1}=\emptyset$.
 assume
\begin{align}
0< \overline r <\underline r < \infty.  \label{rrassum}
\end{align}
Then, recalling Remark \ref{reC:int}, there exist $r'\in {\cal R}_{a_0} $ and $r''\in {\cal R}_{a_1} $ such that
$$
0< \overline r<r'<r''<\underline r< \infty.
$$
By Lemma \ref{lemma:connect}-i), $r'<r''\in {\cal  R}_{a_0}$, and
then $r''\in {\cal R}_{a_0}\cap {\cal  R}_{a_1}=\emptyset$, which is impossible. Therefore, \eqref{rrassum} doe not hold and we conclude $0<\underline r\le \overline r<\infty$.
   We further assume $\underline r=\overline r$. Then
i)-ii) would imply
$\int_0^1  v_{\underline r}(y) Q_0(dy|0) = v_{\underline r}(1)$,
which is impossible since $v_{\underline r}$ is  strictly increasing
on $[0, 1]$ and (A3) holds. This proves iii).

iv) By the definition of $\underline r$ and $\overline r$, it can be shown using \eqref{Vstab} that
$\theta(r)\in (0, 1)$  for $r\in (\underline r, \overline r)$.
By the continuous dependence of the function $v_r(\cdot)$ on $r$ and the method of proving Lemma \ref{lemma:zzb}, we can show  the continuity of $\theta(r)$ on $(0,1)$, and further show $\lim_{r\to \underline r^+}\theta(r)=0$ and $\lim_{r\to \overline{r}^-} \theta(r) =1$. So $\theta(r)$ is continuous on $[\underline r, \overline r]$.
If $\theta(r)$ were not strictly increasing on $[\underline r, \overline r]$,  there would exist
 $\underline r<r_1<r_2<\overline r$ such that
 \begin{align} \theta(r_1)\ge \theta(r_2). \label{s1ges2}
 \end{align}
  If $\theta (r_1)>\theta (r_2)$ in \eqref{s1ges2}, by the continuity of $\theta(r)$,
$\theta(\underline r)=0$, $\theta(\overline r)=1$, and the intermediate value theorem we may find $r'\in (\underline r, r_1)$ such that  $\theta( r_1')=\theta( r_2)$.
Next, we replace $r_1$ by $r_1'$.
Thus if $\theta(r)$ is not strictly increasing,  we may find $r_1<r_2$ from $(\underline r, \overline r)$ such that $\theta (r_1)=
\theta( r_2)\in (0, 1)$, which is a contradiction to Lemma \ref{lemma:the12}.
 This proves iv).
 \qed

\begin{remark}
By Lemmas \ref{lemma:connect} and \ref{lemma:interv},
${\cal R}_{a_1}= (0, \underline r)$ and
${\cal R}_{a_0}= (\overline r, \infty)$.
\end{remark}

\begin{remark}
If (A4) is replaced by $P(\xi\in (0, 1))>0$ without assuming a probability density function $f_\xi$, we still have $C_{R_1}>0$ in the proof of Lemma \ref{lemma:th0}
 and all results in this appendix hold.
\end{remark}

%\begin{remark}
%Based on Lemma \ref{lemma:interv}-iv), we can further apply a
%contradictory argument to show
%$\lim_{r\to \underline{r}^+} \theta(r) =0 $ and
%$\lim_{r\to \overline{r}^-} \theta(r) =1$.
%\end{remark}

%%%%%%%%%%%%%%%%%%%%

\section*{Appendix D: Proof of Theorem \ref{theorem:zmono}}   %6
\label{sec:proofz}

\renewcommand{\theequation}{D.\arabic{equation}}
\setcounter{equation}{0}
\renewcommand{\thetheorem}{D.\arabic{theorem}}
\setcounter{theorem}{0}
\renewcommand{\thelemma}{D.\arabic{lemma}}
\setcounter{lemma}{0}
\renewcommand{\thesubsection}{D.\arabic{subsection}}
\setcounter{subsection}{0}
\renewcommand{\theremark}{D.\arabic{remark}}
\setcounter{remark}{0}

% We combine a vanishing discount

%\subsection{The Long-Run Average of the State}

Let $\{x_t^{i,\theta}, t\ge 0\}$ be the Markov chain generated by  the threshold policy with parameter $0<\theta<1$, where $x_0^{i,\theta}$ is given.
Define $m_t^{i,\theta} = 1-x_t^{i,\theta}$.
By Theorem \ref{theorem:erg},  $\{x_t^{i,\theta}, t\ge 0\}$ and $\{m_t^{i,\theta}, t\ge 0\}$ are  ergodic.

To facilitate further computation, we define an auxiliary  Markov chain $\{Y_t, t\ge 0\}$.
 Let $\xi$ be specified in  (A3) and define $\tilde \xi=1-\xi$.
Let $\{\tilde \xi, \tilde \xi_t, t\ge 1\}$ be i.i.d. random variables.
For $\lambda \in (0, 1)$, define  $\{Y_t, t\ge 0\}$ as follows:
\begin{align}
Y_0=1, \quad Y_t=\tilde\xi_t Y_{t-1} \quad \mbox{for} \quad 1\le t\le
\tau, \label{Yc}
\end{align}
where
$$\tau= \inf\{t| Y_t\le \lambda\}.$$
By (A4),  $P(\tau <\infty)=1$, and
moreover,  $E\tau<\infty$. Set $Y_{\tau+1}=1$ and
the process $\{Y_t, t\ge 0\}$ further evolves from state 1 at time $\tau+1$ as a Markov chain with a stationary transition probability kernel.
Denote
$
S_t =\sum_{i=0}^t Y_i
$ for $t\ge 0$.

\begin{lemma} \label{lemma:YS}
We have
\begin{align}
\lim_{k\to \infty}\frac{1}{k} \sum_{t=0}^{k-1} Y_t= \frac{ES_\tau}
{ 1+E\tau} \qquad{\rm  w.p.1}.  \label{lrun}
\end{align}
\end{lemma}

\proof Take $\lambda=1-\theta$.   Since $\{Y_t, t\ge 0\}$ has the same transition probability kernel as $\{m_t^{i, \theta}, t\ge 0 \}$, it is ergodic, and therefore the left hand side of \eqref{lrun} has a constant limit   w.p.1.
Define $T_0=0$ and $T_n$ as the time for $\{Y_t, t\ge 0\}$ to return to state 1 for the $n$th time. Define  $B_n= \sum_{t= T_{n-1}}^{T_n-1} Y_t$ for $n\ge 1$. We observe that $\{Y_t, t\ge 0\}$ is a regenerative process
(see e.g. \cite{A03,SW93} and \cite[Theorem 4]{AR14})
 with regeneration times $\{T_n, n\ge 1\}$ and that $\{B_n, n\ge 1\}$ is a sequence of i.i.d. random variables. Note that $B_1=S_\tau $ is the sum of $\tau+1$ terms. By the strong  law of large numbers for regenerative processes \cite[pp. 177]{A03},   % we have
%$$
%\lim_{k\to \infty}\frac{1}{k} \sum_{t=0}^{k-1} Y_t = \frac{ES_\tau}{E(\tau+1)}, %\quad \mbox{w.p.1.}
%$$
 the lemma follows. \qed

Define another Markov chain $\{Y_t', t\ge 0\}$ after replacing $\lambda$ in \eqref{Yc} by $\lambda'\in (0, \lambda)$, and $\tau$ by $\tau'=\inf\{t|Y_t'\le \lambda'\}$. The initial state is $Y'_0=1$.  Let
$S'_{\tau'}= \sum_{t=0}^{\tau'} Y_t'$.

\begin{lemma} \label{theorem:SS}
We have
\begin{align}
\frac{ES_\tau}{ 1+E\tau}  \ge \frac{ES'_{\tau'}}{ 1+E\tau'}.
\label{Scc}
\end{align}
\end{lemma}

\proof
Denote $\zeta_k= \prod_{t=1}^k \tilde \xi_t$ for $k\ge 1$.
If $k\ge 2$,  $\{\tau\ge k\}= \{\zeta_{k-1}>\lambda\}$.
For the given $\lambda$, we have
\begin{align}
E\tau= \sum_{k=1}^\infty kP(\tau =k)= \sum_{k=1}^\infty P(\tau \ge k) =1 +\sum_{k=1}^\infty P(\zeta_k>\lambda).
\end{align}

Denote $\alpha=E\tilde \xi$. Then $\alpha \in (0, 1)$ by (A4).
We obtain
\begin{align*}
ES_\tau&=E\sum_{k=1}^\infty  S_k1_{\{\tau=k\}}\\
   % &= 1+\sum_{k=1}^\infty E (\sum_{t=1}^k Y_t)1_{\{\tau =k\}}\\
    &= EY_{0}+EY_1+E\sum_{k=2}^\infty  Y_k 1_{\{\tau\ge k\}} .
\end{align*}

By the independence of $\tilde \xi_k$ and $\zeta_{k-1}1_{\{\zeta_{k-1}>\lambda \}}$ for $k\ge 2$,  it follows that
$$
E(Y_k 1_{\{\tau\ge k\}}) =  E (\tilde \xi_k  \zeta_{k-1}1_{\{ \zeta_{k-1} >\lambda \}})
= \alpha  E (\zeta_{k-1}1_{\{ \zeta_{k-1} >\lambda \}}).
$$
This gives
$$
ES_\tau=
1+\alpha +\alpha\sum_{k=1}^\infty E (\zeta_k 1_{\{\zeta_k>\lambda \}} ).
$$
For $k\ge 1$, denote
$$
p_k= P(\zeta_k>\lambda), \quad r_k=E (\zeta_k 1_{\{\zeta_k>\lambda \}}), \quad
\delta_k= P(\lambda'<\zeta_k\le \lambda),  \quad
\Delta_k=E (\zeta_k 1_{\{ \lambda'< \zeta_k\le \lambda \}} ).
$$
We have
\begin{align}
\frac{ES_\tau}{ 1+E\tau} =
\frac{1+\alpha +\alpha \sum_{k=1}^\infty r_k}
{2+\sum_{k=1}^\infty p_k} \nonumber
\end{align}
and
$$
\frac{ES'_{\tau'}}{ 1+E\tau'}=\frac{1+\alpha +\alpha \sum_{k=1}^\infty (r_k+\Delta_k)}
{2+\sum_{k=1}^\infty (p_k +\delta_k)}.
$$
The inequality \eqref{Scc}
 is equivalent to
\begin{align}
\alpha  \left(\sum_{k=1}^\infty \Delta_k \right) \left(2+\sum_{k=1}^\infty p_k\right) \le \left(1+\alpha +\alpha \sum_{k=1}^\infty r_k \right) \left(\sum_{k=1}^\infty \delta_k\right). \label{eac}
\end{align}
Clearly, $\Delta_k \le \lambda \delta_k $ for $k\ge 1$.  To prove \eqref{eac}, it suffices to show
\begin{align}
\alpha \lambda\left(2+\sum_{k=1}^\infty p_k\right)
\le 1+\alpha +\alpha \sum_{k=1}^\infty r_k. \nonumber
\end{align}
Since $r_k\ge \lambda p_k$  for $k\ge 1$,
we only need to show
$$
2\alpha \lambda \le 1+\alpha,
$$
which follows from
$
0<\alpha<1,\  0<\lambda<1.
$
\qed

Suppose $0<\theta<\theta'<1$. Then there exist two constants $C_\theta, C_{\theta'}$ such that
$$
\lim_{k\to \infty }  \frac{1}{k}\sum_{t=0}^{k-1}x_t^{i,\theta}=C_\theta, \qquad \lim_{k\to \infty }  \frac{1}{k}\sum_{t=0}^{k-1}x_t^{i, \theta'}
= C_{\theta'},  \quad \mbox{w.p.1.}
$$

\begin{lemma} \label{lemma:bb}
We have $C_\theta\le C_{\theta'}$.
\end{lemma}

\proof
The lemma follows from the relation
$x_t^{i,\theta}= 1-m_t^{i, \theta}
$ for $t\ge 0$,    Lemmas  \ref{lemma:YS} and \ref{theorem:SS}.~\qed

\begin{remark} Due to the ergodicity of the Markov chains in  Lemma  \ref{lemma:bb}, the initial   states $x_0^{i, \theta}$ and $x_0^{i, \theta'}$ can be arbitrary.
\end{remark}

{\it Proof of Theorem \ref{theorem:zmono}}.
By the ergodicity of $\{x_t^{i,\theta}, t\ge 0\}$, we have
$
z(\theta_l)= \lim_{k\to \infty }  \frac{1}{k}
\sum_{t=0}^{k-1}x_t^{i, \theta_l}
$ w.p.1. Lemma
\ref{lemma:bb} implies
$z(\theta_1)\le z(\theta_2).$
  \qed

\end{document}